\documentclass[11pt]{article}

\usepackage[T2A]{fontenc}
\usepackage[cp1251]{inputenc}
\usepackage{amsthm}%%%%%%%% theorem definition
\usepackage{amsfonts}
\usepackage{eufrak}
\usepackage{amssymb}
\usepackage{amsmath}
\usepackage{cite}%%%%% citations


\begin{document}
%%%%%%%%%%%%% begin theorem definition %%%%%%%%%%%%%%%%%%
\newtheoremstyle{mytheorem}
  {\topsep}   % ABOVESPACE
  {\topsep}   % BELOWSPACE
  {\itshape}  % BODYFONT
  {}       % INDENT (empty value is the same as 0pt)
  {\bfseries} % HEADFONT
  { }         % HEADPUNCT
  {5pt plus 1pt minus 1pt} % HEADSPACE
  { }          % CUSTOM-HEAD-SPEC
\newtheoremstyle{myremark}
  {\topsep}   % ABOVESPACE
  {\topsep}   % BELOWSPACE
  {\upshape}  % BODYFONT
  {}       % INDENT (empty value is the same as 0pt)
  {\bfseries} % HEADFONT
  { }         % HEADPUNCT
  {5pt plus 0pt minus 1pt} % HEADSPACE
  {}          % CUSTOM-HEAD-SPEC\cite{}
\theoremstyle{mytheorem}
\newtheorem{theorem}{Theorem}[section]
 \newtheorem{theorema}{Theorem}
 \newtheorem*{heyde1*}{Theorem A}
 \newtheorem*{heyde2*}{Theorem B}
 \newtheorem*{heyde3*}{Theorem C}
 \newtheorem*{heyde4*}{Theorem D}
 \newtheorem{proposition}{Proposition}[section]
 \newtheorem{lemma}[theorem]{Lemma}
\newtheorem{corollary}[theorem]{Corollary}
\newtheorem{definition}[theorem]{Definition}
\theoremstyle{myremark}
\newtheorem{remark}[theorem]{Remark}
%%%%%%%%%%%%%%%%%%%%% end theorem definition %%%%%%%%%%%%%%%%%%

\noindent This work with some changes and additions was accepted \\ for publication in the journal "Results in Mathematics"

\vskip 1 cm

\centerline{\bf Generalization of the Heyde  theorem to finite Abelian groups  }

\centerline{\bf   and groups  of the form $\mathbb{R}\times G$,  
where $G$ is a finite Abelian group}

\bigskip

\centerline{\textbf{G.M. Feldman}} 

\bigskip

\noindent{According to the well-known   Heyde theorem the  Gaussian distribution  on the real line is characterized by the symmetry of the conditional  distribution of one linear form of independent random variables given another. We study  analogues of this theorem for some locally compact    Abelian groups. We consider linear forms of two independent random variables with values in a locally compact    Abelian group  $X$. We assume that the characteristic functions of these independent random variables do not vanish. Unlike most previous works, we do not impose any restrictions on coefficients of the linear forms. They are arbitrary topological automorphisms of $X$. }

\bigskip

\noindent {\bf Mathematics Subject Classification.}   43A25,  43A35, 60B15, 62E10.

\bigskip

\noindent{\bf Keywords.} Heyde's theorem,   topological
  automorphism,
locally compact Abelian group.

\bigskip

\section{Introduction}

According to the well-known   Heyde theorem  the   Gaussian distribution  on the real line is characterized by the symmetry of the conditional  distribution of one linear form of independent random variables given another   (\!\!\cite[Theorem 13.4.1]{KaLiRa}).
For two independent random variables   the Heyde theorem  can be formulated as follows.
\begin{heyde1*}  Let $\xi_1$ and $\xi_2$ be real
independent random variables with distributions
$\mu_1$ and $\mu_2$, and let
$a\ne 0$ and $a\ne -1$.  Assume that the conditional  distribution of the linear form $L_2 = \xi_1 +
a\xi_2$ given $L_1 = \xi_1 + \xi_2$ is symmetric.
Then  $\mu_j$  are Gaussian distributions.
\end{heyde1*}
%Note that we also consider degenerate distributions to be Gaussian.
Many studies have been devoted to    analogues of   the Heyde theorem for different classes of locally compact Abelian groups (see e.g. \cite{Fe2, MiFe1, Fe4, Fe20bb, My2, Fe6, FeTVP1, FeTVP, {M2013}, {M2020}, F_solenoid, POTA}). The simplest class of locally compact Abelian groups is the finite Abelian groups. For these groups shifts of the Haar distributions on subgroups   play the role of the Gaussian distributions. Taking this into account, the following statement can be viewed as an analogue of  the Heyde theorem for finite Abelian groups.
%two independent random variables taking values in a  finite Abelian group.
\begin{heyde2*}  {\rm(\!\!\cite[Theorem 1]{Fe2}, see also \cite[Theorem 17.1]{Fe5})} Let $X$ be a finite Abelian group containing no elements of order  $2$.
Let  $\alpha$ be an automorphism of the group $X$, and let $I$ be the identity automorphism. Assume that  $I\pm \alpha$ are also automorphisms of   $X$.
Let $\xi_1$ and $\xi_2$ be
independent random variables with values in   $X$ and distributions
$\mu_1$ and $\mu_2$. Assume that the conditional  distribution of the linear form $L_2 = \xi_1 + \alpha\xi_2$ given $L_1 = \xi_1 + \xi_2$ is symmetric.
Then  $\mu_j$ are shifts of the Haar distribution on a subgroup of  $X$.
\end{heyde2*}
In \S 2 we prove another analogue of   the Heyde theorem for finite Abelian groups.   In contrast to Theorem B
we do not impose any restrictions on an  automorphism $\alpha$ of the group  $X$, but assume that the  characteristic functions of   $\xi_j$ do not vanish. Based on the results of \S 2, we prove  in \S 3 and \S 4 under these assumptions    analogues of    the Heyde theorem for  some other locally compact Abelian groups, namely for   groups of the form $\mathbb{R}\times G$ and $\Sigma_\text{\boldmath $a$}\times G$, where  $G$ is a finite Abelian group containing no elements of order  $2$, and $\Sigma_\text{\boldmath $a$}$ is the  \text{\boldmath $a$}-adic solenoid with \text{\boldmath $a$}=$(2, 3, 4,\dots)$. In the course of the proof  we use some concepts and results from abstract harmonic analysis (see e.g. \cite{Hewitt-Ross}), and also some results on the theory of entire functions, in particular the Hadamard theorem on the representation of an entire function of finite order.

Let $X$ be a second countable locally compact Abelian group. We consider only such groups.   Denote by $Y$ the character
group of the group $X$. For  $x \in X$ denote by  $(x,y)$ the value of a character $y \in Y$ at the element  $x$. For a closed subgroup $L$ of the group $Y$  denote by
 $A(X, L) = \{x \in X: (x, y) =1$ \mbox{for all } $y\in
L\}$
its annihilator. Denote by ${\rm Aut}(X)$ the group
of topological automorphisms of $X$, and by  $I$ the identity automorphism of a group.   Let $\alpha\in{\rm Aut}(X)$ and $K$ be a closed subgroup of   $X$ such that $\alpha(K)=K$, i.e. the restriction  of  $\alpha$ to $K$ is a topological automorphism of the group   $K$. Denote by   $\alpha_{K}$ this restriction.
Let
 $\alpha:X\rightarrow X$ be a topological automorphism of    $X$. The adjoint automorphism
$\tilde\alpha: Y\rightarrow Y$
is defined by the formula $(\alpha x,
y)=(x, \tilde\alpha y)$ for all $x\in X$, $y\in
Y$.  Note that $\alpha\in{\rm Aut}(X)$ if and only if $\tilde\alpha\in{\rm Aut}(Y)$.  Set $X^{(2)}=\{2x: x\in X\}$.   Denote by $\mathbb{R}$ the additive group of real numbers, by $\mathbb{Z}(n)=\{0, 1, \dots, n-1\}$   the group of residue classes modulo $n$, and by $\mathbb{C}$   the   complex plane.

Let $\mu$ be a probability distribution  on   $X$.
Denote by
$$
\hat\mu(y) =
\int_{X}(x, y)d \mu(x), \ \ y\in Y,$$
the characteristic function (Fourier transform) of
the distribution  $\mu$, and by $\sigma(\mu)$ the support of $\mu$. Define the distribution $\bar \mu $ by the formula
 $\bar \mu(B) = \mu(-B)$ for any Borel  subset $B$ of $X$.
Then $\hat{\bar{\mu}}(y)=\overline{\hat\mu(y)}$. For $x\in X$ denote by $E_x$  the degenerate distribution
 concentrated at the element $x$. We need the following definition (see \!\!\cite[Chapter IV, \S 6]{Pa}).

\begin{definition} \label{d1} A distribution  $\gamma$ on a locally compact Abelian group  $X$ is said to be Gaussian
 if its characteristic function is represented in the form
$$
\hat\gamma(y)= (x,y)\exp\{-\varphi(y)\}, \ \  y\in Y,
$$
where $x \in X$, and $\varphi(y)$ is a continuous nonnegative function on the group   $Y$ satisfying the equation
$$
\varphi(u + v) + \varphi(u
- v) = 2[\varphi(u) + \varphi(v)], \ \ u,  v \in
Y.
$$
\end{definition}
Let us make the following general remark. Let $X$ be a locally compact Abelian group, $\xi_1$ and $\xi_2$ be
independent random variables with values in   $X$ and     $\alpha_j, \beta_j\in {\rm Aut}(X)$.   Assume that  the conditional  distribution of the linear form $L_2 = \beta_1\xi_1 + \beta_2\xi_2$ given $L_1 = \alpha_1\xi_1 + \alpha_2\xi_2$ is symmetric. If we are interested in describing of the distributions of $\xi_j$ we can suppose without loss of generality that $L_1 = \xi_1 + \xi_2$ and $L_2 = \xi_1 + \alpha\xi_2$, where $\alpha\in {\rm Aut}(X)$.

\section{ Generalisation of the Heyde theorem to finite Abelian groups}

In this section we   prove the following analogue of   the Heyde theorem for finite Abelian groups.
\begin{theorem}\label{th1} Let $X$ be a finite Abelian group containing no elements of order  $2$, and let  $\alpha$ be an automorphism of the group   $X$. Set $K={\rm Ker}(I+\alpha)$. Let $\xi_1$ and $\xi_2$ be
independent random variables with values in   $X$ and distributions
$\mu_1$ and $\mu_2$ with nonvanishing characteristic functions. Assume that the conditional  distribution of the linear form $L_2 = \xi_1 + \alpha\xi_2$ given $L_1 = \xi_1 + \xi_2$ is symmetric.
Then  $\mu_j=\omega*E_{x_j}$, where $\omega$ is   a distribution   supported in $K$, $x_j\in X$, $j=1, 2$.
%Moreover, if  $\eta_j$ are independent identically distributed random variables with values in   $X$  and distribution  $\omega$, then    the conditional  distribution of the linear form
%$N_2=\eta_1 + \alpha\eta_2$ given $N_1=\eta_1 + \eta_2$ is symmetric.
\end{theorem}
To prove Theorem \ref{th1} we need the following lemmas.
\begin{lemma}\label{lem1}{\rm(\!\!\cite[Lemma 16.1]{Fe5}).} Let $X$ be a locally compact Abelian group,  and let $\alpha$ be a topological automorphism  of the group $X$.
Let
$\xi_1$ and  $\xi_2$  be independent random variables with values in
$X$  and distributions $\mu_1$ and $\mu_2$.   The conditional distribution of the linear form $L_2 = \xi_1 + \alpha\xi_2$ given $L_1 = \xi_1 + \xi_2$ is symmetric if and only
 if the characteristic functions
 $\hat\mu_j(y)$ satisfy the equation
\begin{equation}\label{11.04.1}
\hat\mu_1(u+v )\hat\mu_2(u+\tilde\alpha v )=
\hat\mu_1(u-v )\hat\mu_2(u-\tilde\alpha v), \ \ u, v \in Y.
\end{equation}
\end{lemma}
It is convenient for us to formulate as a  lemma  the following   well-known statement  (see e.g.   \cite[Proposition 2.13]{Fe5}).
\begin{lemma}\label{lem2}  Let $X$ be a locally compact Abelian group, and let $\mu$ be a distribution on $X$. The sets $E=\{y\in Y:  \hat\mu(y)=1\}$ and $B=\{y\in Y:  |\hat\mu(y)|=1\}$ are closed subgroups of the group $Y$,    $\mu$ is supported in the annihilator  $A(X,E)$, and there is an element $x\in X$ such that $\hat\mu(y)=(x, y)$ for all $y\in B$.
\end{lemma}
\begin{lemma}\label{lem9} {\rm(\!\!\cite[(24.22)]{Hewitt-Ross})} Let $X$ be a locally compact Abelian group. The subgroup   $Y^{(2)}$ is dense in the group $Y$ if and only if the group $X$ contains no elements of order  $2$.
\end{lemma}
\begin{lemma}\label{lem11} Let $X$ be a locally compact Abelian group, and
 let $H$ be a closed subgroup of the group  $Y$ such that the subgroup   $H^{(2)}$ is dense in   $H$.   Put   $K=A(X, H)$.
 Let  $\alpha$ be a topological automorphism of the group   $X$ such that $\tilde\alpha(H)=H$.
Let
$\xi_1$ and  $\xi_2$  be independent random variables with values in
 $X$  and distributions $\mu_1$ and $\mu_2$ such that
 \begin{equation}\label{19.04.1}
|\hat\mu_1(y)|=|\hat\mu_2(y)|=1, \ \  y\in H.
\end{equation}
Assume that the conditional  distribution of the linear form $L_2 = \xi_1 + \alpha\xi_2$ given $L_1 = \xi_1 + \xi_2$ is symmetric. Then there are some shifts $\lambda_j$ of the distributions $\mu_j$  such that $\lambda_j$ are supported in      $K$, and if  $\eta_j$ are independent random variables with values in   $X$  and distributions $\lambda_j$, then    the conditional  distribution of the linear form
$N_2=\eta_1 + \alpha\eta_2$ given $N_1=\eta_1 + \eta_2$ is symmetric.
\end{lemma}
{\rmfamily{\bfseries{\itshape Proof}}}   By Lemma \ref{lem1},
the characteristic functions $\hat\mu_j(y)$   satisfy equation (\ref{11.04.1}). By Lemma \ref{lem2},  it follows from (\ref{19.04.1}) that there are elements $p_j\in X$ such that
\begin{equation}\label{11.04.3}
\hat\mu_j(y)=(p_j, y), \ \ y\in H, \ \ j=1, 2.
\end{equation}
Taking into account that $\tilde\alpha(H)=H$, we can  consider the restriction  of equation (\ref{11.04.1}) to $H$ and use  (\ref{11.04.3}).  We get
\begin{equation}\label{11.04.4}
(p_1, u+v)(p_2, u+\tilde\alpha v)=(p_1, u-v)(p_2, u-\tilde\alpha v), \ \ u, v\in H.
\end{equation}
It follows from (\ref{11.04.4}) that
\begin{equation}\label{22.08.1}
2(p_1+\alpha p_2)\in K.
\end{equation}
Note that the character group of the factor-group $X/K$ is topologically isomorphic to the annihilator $A(Y, K)$, and $A(Y, K)=H$. Therefore, by Lemma  \ref{lem9},
if  the subgroup   $H^{(2)}$ is dense in   $H$, then  the factor-group $X/K$ contains no elements of order 2. Hence, the subgroup $K$ has the property: if $2x\in K$, then $x\in K$. For this reason (\ref{22.08.1}) implies that $p_1+\alpha p_2\in K$. Put $x_1=-\alpha p_2$, $x_2=p_2$. Then $x_1+\alpha x_2=0$, and for this reason the characteristic functions $(x_j, y)$ satisfy equation (\ref{11.04.1}). Hence,   the characteristic functions of the distributions  $\lambda_j=\mu_j*E_{-x_j}$ also satisfy equation  (\ref{11.04.1}). Let  $\eta_j$   be independent random variables with values in
the group $X$  and distributions   $\lambda_j$. By Lemma \ref{lem1}, the conditional distribution of the linear form $N_2=\eta_1 + \alpha\eta_2$ given $N_1=\eta_1 + \eta_2$ is symmetric. Taking into account that $-x_1=-p_1+k$, where $k\in K$, $-x_2=-p_2$, and $K=A(X, H)$, we find from (\ref{11.04.3}) that
$\hat\lambda_j(y)=1$ for all $y\in H$, $j=1, 2$.
By Lemma \ref{lem2}, this implies that  $\sigma(\lambda_j)\subset K$, $j=1, 2$. $\hfill\Box$
\begin{lemma}\label{lem3} {\rm(\!\!\cite[Corollary 1]{FeTVP1}).} Let $X$ be a locally compact Abelian group containing no elements of order  $2$.
Let
$\xi_1$ and  $\xi_2$  be independent random variables with values in
$X$  and distributions $\mu_1$ and $\mu_2$.   The conditional distribution of the linear form  $L_2 = \xi_1 -\xi_2$ given
  $L_1 = \xi_1 + \xi_2$ is symmetric if and only
 if  $\mu_1=\mu_2$.
\end{lemma}

\medskip

{\rmfamily{\bfseries{\itshape Proof of Theorem \ref{th1}}}} We note that the group $Y$ is isomorphic to the group $X$, and hence,   the group $Y$ also contains no elements of order   2.

1. First we show     that   the proof of the theorem  is reduced to the case  when $(I-\alpha)\in {\rm Aut}(X)$.
By Lemma \ref{lem1},
the characteristic functions $\hat\mu_j(y)$   satisfy equation (\ref{11.04.1}).
Set $L={\rm Ker}(I-\tilde\alpha)$ and assume that $L\ne\{0\}$. It is obvious that $\tilde\alpha(L)=L$. Since $\tilde\alpha y=y$ for all $y\in L$,   the restriction of equation   (\ref{11.04.1}) to the subgroup $L$ takes the form
\begin{equation}\label{11.04.2}
\hat\mu_1(u+v)\hat\mu_2(u+v)=\hat\mu_1(u-v)\hat\mu_2(u-v), \ \ u, v\in L.
\end{equation}
Substituting  $u=v=y$ in (\ref{11.04.2}), we obtain $\hat\mu_1(2y)\hat\mu_2(2y)=1$ for all $y\in L$. Since the subgroup  $L$  contains no elements of order  $2$, we have ${L^{(2)}}=L$. Thus,    $\hat\mu_1(y)\hat\mu_2(y)=1$, and hence, $|\hat\mu_1(y)|=|\hat\mu_2(y)|=1$ for all $y\in L$. Put $G=A(X, L)$.
By Lemma \ref{lem11},  we can  replace  the distributions  $\mu_j$ by their shifts  $\lambda_j$   in such a way that   $\sigma(\lambda_j)\subset G$, and   if   $\eta_j$
are independent random variables with values in
the group $X$  and distributions $\lambda_j$, then the conditional distribution of the linear form $N_2=\eta_1 + \alpha\eta_2$ given $N_1=\eta_1 + \eta_2$  is symmetric.
It follows from  $\tilde\alpha(L)=L$   that $\alpha(G)=G$, i.e. the restriction  of  $\alpha$ to $G$ is an   automorphism of the group   $G$.   Since $\sigma(\lambda_j)\subset G$, we can consider  $\eta_j$ as independent random variables with values in
the group $G$. In so doing,    the conditional distribution of the linear form  $N_2=\eta_1 + \alpha_{G}\eta_2$ given $N_1=\eta_1 + \eta_2$  is symmetric. Note that $G$ is a proper subgroup of the group $X$ because  $L\ne\{0\}$. Obviously, ${\rm Ker}(I+\alpha_{G})\subset{\rm Ker}(I+\alpha)$. Thus, we may prove the theorem supposing that the independent random variables take values not in the group  $X$ but in
the group  $G$.    Furthermore,  if ${\rm Ker}(I-\tilde\alpha_{G})\ne\{0\}$, then we  repeat this argument. Since  $G$ is a finite group, in a  number of steps, each time replacing the distributions  with their shifts,  we obtain a subgroup $F$ of the group $X$, an automorphism $\alpha_F$ of the group $F$ satisfying the condition     ${\rm Ker}(I-\tilde\alpha_F)= \{0\}$, and independent random variables $\zeta_j$ with values in the group  $F$ such that the conditional distribution of the linear form  $M_2  = \zeta_1 + \alpha_F\zeta_2$ given  $M_1 = \zeta_1 + \zeta_2$ is symmetric. Furthermore,  the distributions of the random variables $\zeta_j$ are shifts of  $\mu_j$. Since  $F$ is a finite group,   ${\rm Ker}(I-\tilde\alpha_F)= \{0\}$ implies that $(I-\alpha_F)\in {\rm Aut}(F)$. Thus, we  can suppose from the beginning that   $(I-\alpha)\in {\rm Aut}(X)$.

Note that in this argument we did not use the fact that the characteristic functions $\hat\mu_j(y)$ do not vanish.

2. Put $H=(I+\tilde\alpha)Y$ and prove that (\ref{19.04.1}) is satisfied.  Set $\nu_j=\mu_j*\bar\mu_j$. Then $\hat\nu_j(y)=|\hat\mu_j(y)|^2>0$  for all $y\in Y$, $j=1, 2$. By Lemma \ref{lem1},
the characteristic functions $\hat\mu_j(y)$   satisfy equation (\ref{11.04.1}). Obviously, the characteristic functions   $\hat\nu_j(y)$ also satisfy equation   (\ref{11.04.1}). Set $f(y)=\hat\nu_1(y)$, $g(y)=\hat\nu_2(y)$ and rewrite   equation (\ref{11.04.1}) for the functions $\hat\nu_j(y)$ in these notation
\begin{equation}\label{11.04.6}
f(u+v)g(u+\tilde\alpha v)=f(u-v)g(u-\tilde\alpha v), \ \ u, v\in Y.
\end{equation}
Substituting  $u=v=y$ in (\ref{11.04.6}), we obtain
\begin{equation}\label{11.04.7}
f(2y)g((I+\tilde\alpha) y)=g((I-\tilde\alpha) y), \ \ y\in Y.
\end{equation}
Since $(I-\alpha)\in {\rm Aut}(X)$, we have $(I-\tilde\alpha)\in {\rm Aut}(Y)$ and   find from    (\ref{11.04.7}) that
\begin{equation}\label{11.04.8}
g(y)=f(2(I-\tilde\alpha)^{-1} y)g((I+\tilde\alpha)(I-\tilde\alpha)^{-1} y), \ \ y\in Y.
\end{equation}
This implies the inequality
\begin{equation}\label{11.04.9}
g(y)\le f(2(I-\tilde\alpha)^{-1} y), \ \ y\in Y.
\end{equation}
Substituting $u=\tilde\alpha y$, $v=y$ in (\ref{11.04.6}), we obtain
\begin{equation}\label{11.04.10}
f((I+\tilde\alpha) y)g(2\tilde\alpha y)=f(-(I-\tilde\alpha) y), \ \ y\in Y.
\end{equation}
It follows from this that
\begin{equation}\label{11.04.11}
f(y)\le g(-2\tilde\alpha(I-\tilde\alpha)^{-1} y), \ \ y\in Y.
\end{equation}
From (\ref{11.04.9}) and (\ref{11.04.11}) we find
\begin{equation}\label{11.04.12}
g(y)\le f(2(I-\tilde\alpha)^{-1} y)\le g(-4\tilde\alpha(I-\tilde\alpha)^{-2} y), \ \ y\in Y.
\end{equation}
 Set $\kappa=-4\tilde\alpha(I-\tilde\alpha)^{-2}$. Since   multiplication by 2 is an automorphism  of the group   $Y$, we have $\kappa\in {\rm Aut}(Y)$. The group ${\rm Aut}(Y)$ is finite because the group $Y$ is finite. Denote by   $m$ the order of the element    $\kappa$ of the group ${\rm Aut}(Y)$.  We find from  (\ref{11.04.12})  the inequalities
$$
g(y)\le g(\kappa y)\le \dots \le g(\kappa^{m-1} y)\le g(\kappa^{m} y)=g(y), \ \ y\in Y.
$$
Hence,
\begin{equation}\label{11.04.13}
g(y)=g(\kappa y), \ \ y\in Y,
\end{equation}
and the  equality
\begin{equation}\label{11.04.14}
g(y)=f(2(I-\tilde\alpha)^{-1} y), \ \ y\in Y,
\end{equation}
follows from (\ref{11.04.12}) and (\ref{11.04.13}). Taking into account that the characteristic functions $f(y)$ and $g(y)$ do not vanish,   and the fact that  $(I-\tilde\alpha)\in{\rm Aut}(Y)$, we find from   (\ref{11.04.8}) and (\ref{11.04.14}) that $g(y)=1$ for all $y\in H$.

Prove an analogous statement for the function   $f(y)$. It follows from (\ref{11.04.9}) and (\ref{11.04.11}) that
\begin{equation}\label{11.04.17}
f(y)\le g(-2\tilde\alpha(I-\tilde\alpha)^{-1} y)\le f(\kappa y), \ \ y\in Y.
\end{equation}
We get from (\ref{11.04.17})  the inequalities
$$
f(y)\le f(\kappa y)\le \dots \le f(\kappa^{m-1} y)\le f(\kappa^{m} y)=f(y), \ \ y\in Y.
$$
Hence,
\begin{equation}\label{11.04.18}
f(y)=f(\kappa y), \ \ y\in Y.
\end{equation}
We obtain from (\ref{11.04.17}) and (\ref{11.04.18})  that
\begin{equation}\label{11.04.15}
f(y)=g(-2\tilde\alpha(I-\tilde\alpha)^{-1} y), \ \ y\in Y.
\end{equation}
Since   $(I-\tilde\alpha)\in {\rm Aut}(Y)$, we find from (\ref{11.04.10}) that
\begin{equation}\label{11.04.16}
f(y)=f(-(I+\tilde\alpha)(I-\tilde\alpha)^{-1} y)g(-2\tilde\alpha (I-\tilde\alpha)^{-1} y), \ \ y\in Y.
\end{equation}
Taking into account that the characteristic functions $f(y)$ and $g(y)$ do not vanish,   and the fact that  $(I-\tilde\alpha)\in{\rm Aut}(Y)$, we find from   (\ref{11.04.15}) and (\ref{11.04.16}) that $f(y)=1$ for all $y\in H$. Thus, we proved that   $f(y)=g(y)=1$ for all $y\in H$, and this implies that
 (\ref{19.04.1}) is satisfied.

3. It follows from $H=(I+\tilde\alpha)Y$ that $K=A(X, H)$.   Since the subgroup  $H$  contains no elements of order  $2$, we have  ${H^{(2)}}=H$. It is obvious that $\tilde\alpha(H)=H$. By Lemma \ref{lem11}, we can replace the distributions $\mu_j$ by  their shifts $\lambda_j$ in such a way   that the distributions   $\lambda_j$ are supported in     the subgroup $K$, and if  $\eta_j$ are independent random variables with values in the group $X$  and distributions $\lambda_j$, then    the conditional  distribution of the linear form
$N_2=\eta_1 + \alpha\eta_2$ given $N_1=\eta_1 + \eta_2$ is symmetric.  In view of  $K={\rm Ker}(I+\alpha)$, the restriction of the automorphism $\alpha$  to the subgroup $K$ coincides with $-I$.
Since $\sigma(\lambda_j)\subset K$, we can consider $\eta_j$ as independent random variables with values in the group   $K$. In so doing, the conditional distribution of the linear form  $N_2=\eta_1 -\eta_2$ given $N_1=\eta_1 + \eta_2$  is symmetric. Applying Lemma \ref{lem3}  to the group   $K$, we get $\lambda_1=\lambda_2=\omega$. This implies the statement of the theorem. $\hfill\Box$
\begin{remark}\label{r1} Let $X$ be a finite Abelian group, and let  $\alpha$ be an automorphism of  the group  $X$. Set $K={\rm Ker}(I+\alpha)$. Let $\omega$ be a  distribution on  $X$ supported in the subgroup $K$, and let $x_1$ and $x_2$ be elements of the group  $X$ such that $x_1+\alpha x_2=0$.  Put $\mu_j=\omega*E_{x_j}$, $j=1, 2$. Let $\xi_j$   be
independent random variables with values in the group $X$ and distributions
$\mu_j$. It is easy to see that   the characteristic functions $\hat\mu_j(y)$ satisfy equation (\ref{11.04.1}). By Lemma  \ref{lem1},  the conditional  distribution of the linear form $L_2 = \xi_1 + \alpha\xi_2$ given $L_1 = \xi_1 + \xi_2$ is symmetric.
From what has been said it follows that    Theorem \ref{th1} can not be strengthened by narrowing the class of distributions which are characterized by the symmetry of the conditional distribution of one linear form given another.
\end{remark}
We give now another proof of Theorem \ref{th1}. It is based on the finite difference method. This method is often used to solve functional equations that arise in characterization problems in  mathematical statistics (see e.g. \cite[\S 10 and  \S 16]{Fe5}).

Let $P(y)$ be a  function on the group   $Y$ and $h\in Y$.
Denote by $\Delta_h$ the finite difference operator
$$\Delta_h P(y)=P(y+h)-P(y), \ \ y\in Y.$$

We need the following statement  which is a special case of a general theorem on continuous polynomials on locally compact Abelian groups (see e.g. \cite[Proposition 5.7]{Fe5}).
\begin{lemma}\label{lem6}  Let $Y$ be a finite   Abelian group,    and let a function $P(y)$ satisfy the equation
$$
 \Delta_{h}^{n}P(y)=0, \ \ y,h \in Y,
$$
for some natural $n$.  Then $P(y)=const$.
\end{lemma}
{\rmfamily{\bfseries{\itshape Second proof of Theorem \ref{th1}}}} We retain the notation used in the proof of Theorem   \ref{th1}. As established in item   1 of the proof of Theorem \ref{th1}, we can suppose that   $(I-\alpha)\in {\rm Aut}(X)$. Then  $(I-\tilde\alpha)\in {\rm Aut}(Y)$. We shall prove that (\ref{19.04.1}) is satisfied.
  This is the main part of the proof of Theorem  \ref{th1}.

Set $P(y)=\log f(y)$, $Q(y)=\log g(y)$. We find from equation (\ref{11.04.6}) that the functions $P(y)$ and $Q(y)$ satisfy the equation
\begin{equation}\label{17.04.1}
P(u+v)+Q(u+\tilde\alpha v)-P(u-v)-Q(u-\tilde\alpha v)=0, \ \ u, v\in Y.
\end{equation}

Equation (\ref{17.04.1})  arises in the study of    the Heyde theorem on various locally compact Abelian groups with certain restrictions on  $\alpha$ (see e.g. \cite{Fe6}, \cite{FeTVP1}, \cite{F_solenoid}). For completeness, we give here its solution.

Let $k_1$ be an arbitrary element of the group $Y$. Put $h_1
=\tilde\alpha k_1$.   Replacing $u$ by $u+h_1$  and $v$ by $v+k_1$ in (\ref{17.04.1})     and  subtracting   equation
(\ref{17.04.1}) from the obtained equation, we find
\begin{equation}\label{17.04.2}
\Delta_{l_{11}}P(u+v) +
\Delta_{l_{12}}Q(u+\tilde\alpha v) -
\Delta_{l_{13}}P(u-v)=0, \ \ u, v \in
Y,
\end{equation}
where $l_{11}= (I+\tilde\alpha)k_1$, $l_{12}=2 \tilde\alpha k_1$, $l_{13}= (\tilde\alpha-I)k_1.$
 Let $k_2$ be an arbitrary element of the group $Y$.   Replacing  $u$ by $u+k_2$  and $v$ by $v+k_2$ in (\ref{17.04.2}) and  subtracting   equation
(\ref{17.04.2}) from the obtained equation, we get
\begin{equation}\label{17.04.3}
\Delta_{l_{21}}\Delta_{l_{11}}P(u+v) +
\Delta_{l_{22}}\Delta_{l_{12}}Q(u+\tilde\alpha v) =
0,
 \ \ u, v \in Y,
\end{equation}
where $l_{21}=2k_2$, $l_{22}=
(I+\tilde\alpha)k_2.$ Let $k_3$ be an arbitrary element of the group $Y$. Put $h_3 =-\tilde\alpha
k_3$.
Replacing $u$ by $u+h_3$  and $v$ by $v+k_3$ in (\ref{17.04.3}) and  subtracting   equation
(\ref{17.04.3}) from the obtained equation, we find
\begin{equation}\label{17.04.4}
\Delta_{l_{31}}\Delta_{l_{21}}\Delta_{l_{11}}P(u+v)
= 0,
 \ \ u, v \in Y,
\end{equation}
where $l_{31}= (I-\tilde\alpha)k_3.$ Substituting $v=0$ in (\ref{17.04.4}), we obtain
\begin{equation}\label{17.04.5}
\Delta_{l_{31}}\Delta_{l_{21}}\Delta_{l_{11}}P(u) = 0,
 \ \ u \in Y.
\end{equation}
Substituting expressions for $l_{11}$, $l_{21}$ and $l_{31}$ into (\ref{17.04.5}) we get
\begin{equation}\label{05.03.1}
\Delta_{(I-\tilde\alpha)k_3}\Delta_{2k_2}\Delta_{(I+\tilde\alpha)k_1}P(u) = 0,
 \ \ u \in Y.
\end{equation}

Since   multiplication by 2 and $I-\tilde\alpha$ are   automorphisms of   the group   $Y$,  and $k_j$ are arbitrary elements of   $Y$, it follows from (\ref{05.03.1})
   that the function $P(y)$ satisfies the equation
$$
\Delta_h^3P(y)=0, \ \ y, h \in H.
$$
Since $P(0)=0$, by Lemma \ref{lem6}, $P(y)=0$ for all $y\in H$. Hence, $f(y)=1$ for all $y\in H$. Reasoning in the same way, we exclude from   equation (\ref{17.04.3}) the function $P(y)$  and get that $Q(y)=0$ for all $y\in H$. Hence, $g(y)=1$ for all $y\in H$. Thus, we proved that   $f(y)=g(y)=1$ for all $y\in H$, and this implies that  (\ref{19.04.1}) is satisfied. The final part of the proof is the same as in item 3 of the proof of Theorem  \ref{th1}.

$\hfill\Box$

\section{ Generalization  of the Heyde  theorem to   groups of the form  $\mathbb{R}\times G$, where $G$ is a finite Abelian group}

Let $X=\mathbb{R}\times G$, where $G$ is a finite Abelian group. Denote by $x=(t, g)$, where $t\in \mathbb{R}$, $g\in G$, elements of the group    $X$.  Let $H$ be the character group of the group    $G$. The   group  $Y$ is topologically isomorphic to the group   $\mathbb{R}\times H$. Denote by $y=(s, h)$, where $s\in \mathbb{R}$, $h\in H$, elements of the group $Y$. Let  $\alpha$ be a topological automorphism of the group
$X$.  It is obvious that $\alpha(\mathbb{R})=\mathbb{R}$ and  $\alpha(G)=G$. Therefore,   $\alpha$ acts on elements of $X$ as follows $\alpha(t, g)=(\alpha_{\mathbb{R}} t, \alpha_{G} g)$, where $t\in \mathbb{R}$, $g\in G$. In so doing,  $\alpha_{\mathbb{R}}$ is multiplication by a nonzero real number $a$.   We will identify $\alpha_{\mathbb{R}}$ with the   number  $a$. The adjoint automorphism   $\tilde\alpha_{\mathbb{R}}$ we shall also   identify with $a$. Thus, $\alpha(t, g)=(a t, \alpha_{G} g)$, and we will write    $\alpha$   in the form $\alpha=(a, \alpha_{G})$.

We will prove the following characterization theorem which  generalizes Theorem \ref{th1}.
\begin{theorem}\label{th2} Let $X=\mathbb{R}\times G$, where $G$ is a finite Abelian group containing no elements of order  $2$.
Let  $\alpha=(a, \alpha_{G})$ be a topological automorphism of the group
  $X$.  Set $K={\rm Ker}(I+\alpha_{G})$.
Let $\xi_1$ and $\xi_2$ be
independent random variables with values in   $X$ and distributions
$\mu_1$ and $\mu_2$ with nonvanishing characteristic functions. Assume that the conditional  distribution of the linear form $L_2 = \xi_1 + \alpha\xi_2$ given $L_1 = \xi_1 + \xi_2$ is symmetric.
If  $a \ne-1$, then $\mu_j=\gamma_j*\omega*E_{x_j}$, where $\gamma_j$ is a Gaussian distribution on   $\mathbb{R}$, $\omega$ is a distribution supported in    $K$, $x_j\in X$, $j=1, 2$.  If  $a=-1$, then $\mu_j=\omega*E_{x_j}$, where  $\omega$ is a distribution supported in  $\mathbb{R}\times K$, $x_j\in X$, $j=1, 2$.
\end{theorem}
To prove Theorem \ref{th2}  we need the following easily verifiable statement  which we formulate as a lemma (see e.g.  \cite[Lemma 6.9]{Fe9}).
\begin{lemma}\label{lem5}
Let $X=\mathbb{R}\times G$, where $G$ is a locally compact Abelian group, and let    $H$ be the character group of $G$. Let
$\mu$ be a distribution on $X$ such that the function $\hat\mu(s, 0)$, $s\in \mathbb{R}$, can be extended to the complex plane $\mathbb{C}$ as an entire function in   $s$. Then  for each fixed   $h\in H$ the function $\hat\mu(s, h)$, $s\in \mathbb{R}$, also can be extended to the complex plane $\mathbb{C}$ as an entire function in     $s$, and the inequality
$$
\max_{s\in \mathbb{C}, \ |s|\le r}|\hat\mu(s, h)|\le \max_{s\in \mathbb{C}, \ |s|\le r}|\hat\mu(s, 0)|
$$
holds.
\end{lemma}
The proof of Theorem \ref{th2} is based on Theorem \ref{th1} and on the following statement  which plays a key role in the proof of Theorem \ref{th2}.
\begin{lemma}\label{lem10}  Let $X=\mathbb{R}\times K$, where $K$ is a locally compact Abelian group. Denote by $(t, k)$, where $t\in \mathbb{R}$, $k\in K$, elements of the group $X$, and by $L$ the character group of the group $K$. Assume that $L^{(2)}=L$.
Let  $\alpha$ be a topological automorphism of the group
$X$ of the form $\alpha(t, k)=(at, -k)$, where $a\ne -1$.
Let $\xi_1$ and $\xi_2$ be
independent random variables with values in   $X$ and distributions
$\mu_1$ and $\mu_2$ with nonvanishing characteristic functions. Assume that the conditional  distribution of the linear form $L_2 = \xi_1 + \alpha\xi_2$ given $L_1 = \xi_1 + \xi_2$ is symmetric.
Then $\mu_j=\gamma_j*\omega$, where $\gamma_j$ is a Gaussian distribution on   $\mathbb{R}$, and $\omega$ is a distribution supported in   $K$,   $j=1, 2$.
\end{lemma}
{\rmfamily   {\bfseries   {\itshape Proof}}}     The group $Y$ is topologically isomorphic to the group   $\mathbb{R}\times L$. Denote by $y=(s, l)$,  where $s\in \mathbb{R}$, $l\in L$, elements of the group $Y$.  By Lemma \ref{lem1}, the characteristic functions   $\hat\mu_j(s, l)$   satisfy equation (\ref{11.04.1})  
which takes the form
$$
\hat\mu_1(s_1+s_2, l_1+l_2)\hat\mu_2(s_1+a s_2, l_1- l_2)$$\begin{equation}\label{14.04.2}=
\hat\mu_1(s_1-s_2, l_1-l_2)\hat\mu_2(s_1-a s_2, l_1+l_2), \ \  s_j\in \mathbb{R}, \ \ l_j\in L.
\end{equation}
Substitute $l_1=l_2=0$ in (\ref{14.04.2}). In view of Lemma \ref{lem1} and Theorem A, the obtained equation implies that
$$
\hat\mu_j(s, 0)=\exp\{-\sigma_js^2+i\beta_j s\}, \ \ s\in \mathbb{R},
$$
where $\sigma_j\ge 0$, $\beta_j\in \mathbb{R}$, $j=1, 2$. Since $\beta_1+a\beta_2=0$, we can replace the distributions    $\mu_j$ by their shifts $\lambda_j=\mu_j*E_{-\beta_j}$ and suppose from the   beginning,   without loss of generality,      that  $\beta_1=\beta_2=0$. We also note that  $\sigma_1+a\sigma_2=0$. This implies that either  $\sigma_1=\sigma_2=0$ or $\sigma_1>0$ and $\sigma_2>0$.

If $\sigma_1=\sigma_2=0$, then $\hat\mu_1(s, 0)=\hat\mu_2(s, 0)=1$ for all $s\in \mathbb{R}$. By Lemma \ref{lem2}, it follows from this that $\sigma(\mu_j)\subset K$, $j=1, 2$. Hence, we may consider $\xi_j$ as independent random variables with values in the group   $K$. If so doing  the conditional distribution of the linear form $L_2=\xi_1 - \xi_2$ given
$L_1=\xi_1 + \xi_2$ is symmetric.    Since $L^{(2)}=L$, by Lemma \ref{lem9}, the group $K$ contains no elements of order   2. The statement of the lemma follows in this case from Lemma \ref{lem1} and Lemma   \ref{lem3}, applied  to the group $K$.

Thus, we will assume that
 \begin{equation}\label{14.04.3}
\hat\mu_j(s, 0)=\exp\{-\sigma_js^2\}, \ \ s\in \mathbb{R},
\end{equation}
where $\sigma_j> 0$, $j=1, 2$.  It follows from  (\ref{14.04.3}) that $a<0$, because if  $a>0$, then, obviously, $\sigma_1=\sigma_2=0$. In view of (\ref{14.04.3}), by Lemma \ref{lem5}, for each fixed   $l\in L$ the functions $\hat\mu_j(s, l)$ can be extended to the complex plane $\mathbb{C}$ as entire functions in   $s$. It is easy to see that equation   (\ref{14.04.2})  holds for all $s_j\in \mathbb{C}$, $l_j\in L$. It also follows from (\ref{14.04.2}) that the functions $\hat\mu_j(s, l)$ do not vanish for all $s\in \mathbb{C}$, $l\in L$.  Indeed, consider   equation (\ref{14.04.2}), assuming that  $s_j\in \mathbb{C}$, $l_j\in L$, and substitute $l_1=l_2=l$ in it.   In view of (\ref{14.04.3}), and taking into account that $L^{(2)}=L$,   rewrite the obtained equation in the form
$$
\hat\mu_1(s_1+s_2, l)\exp\{-\sigma_2(s_1+a s_2)^2\}$$\begin{equation}\label{14.04.4}=
\exp\{-\sigma_1(s_1-s_2)^2\}\hat\mu_2(s_1-a s_2, l), \ \  s_j\in \mathbb{C}, \ \ l\in L.
\end{equation}
Suppose that $\hat\mu_1(s_0, l_0)=0$ for some $s_0\in \mathbb{C}$, $l_0\in L$. Taking into account that $a\ne  -1$, substitute $s_1=s_0a(a+1)^{-1}$, $s_2=s_0(a+1)^{-1}$, $l=l_0$ in   (\ref{14.04.4}). Since $\hat\mu_2(0, l_0)\ne 0$, the right-hand side of the obtained equality is nonzero, and the left-hand side is equal to zero. From the obtained contradiction it follows that the function $\hat\mu_1(s, l)$ does not vanish   for all $s\in \mathbb{C}$, $l\in L$.   Similarly, if $\hat\mu_2(s_0, l_0)=0$ for some $s_0\in \mathbb{C}$, $l_0\in L$, then substituting $s_1=s_0(a+1)^{-1}$, $s_2=-s_0(a+1)^{-1}$, $l=l_0$ in (\ref{14.04.4}),  and taking into account that $\hat\mu_1(0, l_0)\ne 0$, we obtain the contradiction. Hence, the function $\hat\mu_2(s, l)$ also does not vanish   for all $s\in \mathbb{C}$, $l\in L$. It follows from  Lemma \ref{lem5} and  (\ref{14.04.3}) that the functions $\hat\mu_j(s, l)$  for each fixed
$l\in L$ are entire functions in    $s$  of the order at most   2.  By the Hadamard theorem on the representation of an entire function of finite order, taking into account that the functions $\hat\mu_j(s, l)$ do not vanish,  we obtain the representations
\begin{equation}\label{14.04.5}
\hat\mu_j(s, l)=\exp\{a_{lj}s^2+b_{lj}s+c_{lj}\}, \ \ s\in \mathbb{R}, \ \ l\in L, \ \ j=1, 2,
\end{equation}
where $a_{lj}$, $b_{lj}$, $c_{lj}$ are some complex constants. Actually, the representation (\ref{14.04.5}) is valid when $s\in \mathbb{C}$, but we need it only when $s\in \mathbb{R}$. Set
$$
\psi_j(s, l)=  a_{lj}s^2+b_{lj}s+c_{lj}, \ \ s\in \mathbb{R}, \ \ j=1, 2,
$$
and substitute (\ref{14.04.5}) into (\ref{14.04.2}). It follows from the obtained equation that the functions $\psi_j(s, l)$ satisfy the equation
$$
\psi_1(s_1+s_2, l_1+l_2)+\psi_2(s_1+a s_2, l_1- l_2)$$\begin{equation}\label{14.04.6}=
\psi_1(s_1-s_2, l_1-l_2)+\psi_2(s_1-a s_2, l_1+l_2)+2\pi in(l_1, l_2), \ \  s_j\in \mathbb{R}, \ \ l_j\in L,
\end{equation}
where $n(l_1, l_2)$ is an integer. Substituting $l_1=l_2=l$ in (\ref{14.04.6}), we find
$$
a_{2l1}(s_1+s_2)^2+b_{2l1}(s_1+s_2)+c_{2l1}+a_{02}(s_1+as_2)^2
+b_{02}(s_1+as_2)+c_{02} $$$$=a_{01}(s_1-s_2)^2+b_{01}(s_1-s_2)+c_{01}+a_{2l2}(s_1-as_2)^2
+b_{2l2}(s_1-as_2)+c_{2l2}$$\begin{equation}\label{14.04.7}+2\pi in(l, l), \ \ s_j\in \mathbb{R}, \ \ l\in L.
\end{equation}

Equating  the coefficients  of
$s_1^2$ and  $s_2^2$ on each  side  of (\ref{14.04.7}), we obtain
\begin{equation}\label{14.04.8}
a_{2l1}+a_{02}=a_{01}+a_{2l2}, \ \ a_{2l1}+a^2a_{02}=a_{01}+a^2a_{2l2}.
\end{equation}
It follows from (\ref{14.04.3}) that
\begin{equation}\label{23.02.1}
a_{01}=-\sigma_1, \ \ a_{02}=-\sigma_2.
\end{equation}
Taking   into account   that $L^{(2)}=L$  and  $a<0$,   $a\ne -1$,   it follows from  (\ref{14.04.8}) and (\ref{23.02.1}) that
\begin{equation}\label{15.04.1}
a_{l1}= -\sigma_1, \ \ a_{l2}= -\sigma_2, \ \ l\in L.
\end{equation}

Equating  the coefficients  of $s_1$ and $s_2$ on each side  of (\ref{14.04.7}), we find
\begin{equation}\label{14.04.9}
b_{2l1}+b_{02}=b_{01}+b_{2l2}, \ \ b_{2l1}+ab_{02}=-b_{01}-ab_{2l2}.
\end{equation}
It follows from (\ref{14.04.3}) that
\begin{equation}\label{23.02.2}
b_{01}=b_{02}=0.
\end{equation}
Taking   into account   that $L^{(2)}=L$  and  $a<0$,   $a\ne -1$,   it follows from  (\ref{14.04.9}) and (\ref{23.02.2}) that
\begin{equation}\label{15.04.2}
b_{l1}=b_{l2}=0, \ \ l\in L.
\end{equation}

Equating the constant terms on each side  of  (\ref{14.04.7}), we obtain
\begin{equation}\label{14.04.10}
 c_{2l1}+c_{02}=c_{01}+c_{2l2}+2\pi in(l, l), \ \ l\in L.
\end{equation}
In view of (\ref{14.04.3}), we can assume that $c_{01}=c_{02}=0$. Taking this into account  and  the fact that  $L^{(2)}=L$, the equality $\exp\{c_{l1}\}=\exp\{c_{l2}\}$ for all $l\in L$  follows from  (\ref{14.04.10}). We note that (\ref{14.04.5}) implies that $\exp\{c_{lj}\}=\hat\mu_j(0, l)$, $l\in L$,  $j=1, 2$.
Denote by $\omega$ a distribution on the group   $K$ with the characteristic function
\begin{equation}\label{15.04.3}
\hat\omega(l)=\exp\{c_{l1}\}=\exp\{c_{l2}\}, \ \ l\in L.
\end{equation}
Taking into account (\ref{15.04.1}), (\ref{15.04.2}) and (\ref{15.04.3}),   the representation
\begin{equation}\label{14.04.11}
\hat\mu_j(s, l)=\exp\{-\sigma_j s^2\}\hat\omega(l), \ \ s\in \mathbb{R}, \ \ l\in L, \ \ j=1, 2,
\end{equation}
follows from (\ref{14.04.5}).
Denote by $\gamma_j$   Gaussian distributions on  $\mathbb{R}$ with the characteristic functions
$$
\hat\gamma_j(s)=\exp\{-\sigma_js^2\}, \ \ j=1, 2.
$$
It follows from (\ref{14.04.11}) that $\mu_j=\gamma_j*\omega$,  $j=1, 2$, that proves the lemma.

Observe that in the case when $K=\mathbb{Z}(n)$, where $n$ is odd, the lemma was proved in    \cite{F_solenoid}.
$\Box$

\medskip

{\rmfamily   {\bfseries   {\itshape Proof of Theorem \ref{th2}}}} By Lemma \ref{lem1}, the characteristic functions   $\hat\mu_j(s, h)$   satisfy equation (\ref{11.04.1})  which takes the form
$$
\hat\mu_1(s_1+s_2, h_1+h_2)\hat\mu_2(s_1+a s_2, h_1+\tilde\alpha_{G} h_2)$$\begin{equation}\label{13.04.2}=
\hat\mu_1(s_1-s_2, h_1-h_2)\hat\mu_2(s_1-a s_2, h_1-\tilde\alpha_{G} h_2), \ \  s_j\in \mathbb{R}, \ \ h_j\in H.
\end{equation}
 Denote by $\omega_j$   distributions on the group $G$ with the characteristic functions  $\hat\omega_j(h)=\hat\mu_j(0, h)$,  $j=1, 2$. Substituting  $s_1=s_2=0$ in (\ref{13.04.2}), we obtain
\begin{equation}\label{14.04.1}
\hat\omega_1(h_1+h_2)\hat\omega_2(h_1+\tilde\alpha_{G} h_2)=
\hat\omega_1(h_1-h_2)\hat\omega_2(h_1-\tilde\alpha_{G} h_2), \ \    h_j\in H.
\end{equation}
Taking into account    Lemma     \ref{lem1} and equation (\ref{14.04.1}), it follows from Theorem \ref{th1} applying to the group  $G$ and the automorphism $\alpha_G$,  that there are elements
$g_j\in G$ such that   the distribution     $\omega=\omega_1*E_{-g_1}=\omega_2*E_{-g_2}$ is supported in the subgroup $K$.  Moreover,  if  $\zeta_j$ are independent identically distributed random variables with values in the group  $G$ and distribution
    $\omega$, then the conditional distribution of the linear form
$M_2 = \zeta_1 +\alpha_{G}\zeta_2$ given $M_1 = \zeta_1 + \zeta_2$  is symmetric. Consider the distributions  $\lambda_j=\mu_j*E_{-g_j}$.  Let $\eta_j$ be independent random variables with values in the group $X$ and distributions  $\lambda_j$. Then the conditional distribution of the linear form
$N_2=\eta_1 + \alpha\eta_2$ given $N_1=\eta_1 + \eta_2$ is also symmetric. It is obvious that  $\sigma(\lambda_j)\subset\mathbb{R}\times K$. Since  $K={\rm Ker}(I+\alpha_{G})$, the restriction of the automorphism   $\alpha_{G}$  to the subgroup $K$ coincides with $-I$. We can consider  $\eta_j$ as
 independent random variables with values in the subgroup $\mathbb{R}\times K$.
From what has been said  it follows that proving the theorem  we can assume that $X=\mathbb{R}\times K$, where $K$ is a finite Abelian group containing no elements of order 2, and the topological automorphism   $\alpha$   is of the form $\alpha(t, k)=(at, -k)$, where $t\in \mathbb{R}$, $k\in K$.

Let $a\ne - 1$. Denote by $L$ the character group of the group    $K$.    Since the group $K$ contains no elements of order 2, and $L$ is isomorphic to $K$, the group   $L$ also contains no elements of order 2. Hence,  $L^{(2)}=L$. Then the statement of the theorem follows from Lemma  \ref{lem10}.

Let $a=-1$. Then $\alpha=-I$, and by Lemma   \ref{lem3}, we obtain that $\mu_1=\mu_2=\omega$. The theorem is completely  proved. $\hfill\Box$

Note that if in Theorem \ref{th2} $a\ne -1$, then ${\rm Ker}(I+\alpha)=K$, and if $a=-1$, then ${\rm Ker}(I+\alpha)=\mathbb{R}\times K$. Thus, Theorem \ref{th2} implies the following statement.
\begin{corollary}\label{co1} Let $X=\mathbb{R}\times G$, where $G$ is a finite Abelian group containing no elements of order  $2$.
Let  $\alpha$ be a topological automorphism of the group
  $X$.
Let $\xi_1$ and $\xi_2$ be
independent random variables with values in   $X$ and distributions
$\mu_1$ and $\mu_2$ with nonvanishing characteristic functions. If the conditional  distribution of the linear form $L_2 = \xi_1 + \alpha\xi_2$ given $L_1 = \xi_1 + \xi_2$ is symmetric, then
$\mu_j=\gamma_j*\omega$, where $\gamma_j$ is a Gaussian distribution on   $X$, $\omega$ is a distribution supported in      ${\rm Ker}(I+\alpha)$,  $j=1, 2$.
\end{corollary}
\begin{remark}\label{nr1}
 Let $X=\mathbb{R}\times G$, where $G$ is a finite Abelian group. Let  $\alpha=(a, \alpha_{G})$ be a topological automorphism of the group
  $X$. Put $K={\rm Ker}(I+\alpha_G)$.

Assume that $a \ne-1$. Let $\gamma_j$ be   Gaussian distributions on $\mathbb{R}$ with the characteristic functions $\hat\gamma_j(s)=\exp\{-\sigma_js^2\}$, $s\in \mathbb{R}$, where $\sigma_1+a\sigma_2=0$, let $\omega$ be a  distribution on  $X$ supported in the subgroup $K$, and let $x_j$  be elements of the group  $X$ such that $x_1+\alpha x_2=0$. Put $\mu_j=\gamma_j*\omega*E_{x_j}$, $j=1, 2$. Let $\xi_j$  be
independent random variables with values in the group $X$ and distributions
$\mu_j$.  It is easy to see that   the characteristic functions $\hat\mu_j(y)$ satisfy equation (\ref{11.04.1}). By Lemma  \ref{lem1},  the conditional  distribution of the linear form $L_2 = \xi_1 + \alpha\xi_2$ given $L_1 = \xi_1 + \xi_2$ is symmetric.

Assume that   $a=-1$. Then  we argue in the same way as in Remark \ref{r1} and construct independent random variables $\xi_j$   with values in the group $X$ and distributions
$\mu_j=\omega*E_{x_j}$, where   $\omega$ is a distribution supported in    $\mathbb{R}\times K$, $x_j\in X$, $j=1, 2$, such that the conditional  distribution of the linear form $L_2 = \xi_1 + \alpha\xi_2$ given $L_1 = \xi_1 + \xi_2$ is symmetric.

From what has been said it follows that Theorem \ref{th2} can not be strengthened   by narrowing the class of distributions which are characterized by the symmetry of the conditional distribution of one linear form given another.
\end{remark}
\begin{remark}\label{r2}
Let $X=\mathbb{R}\times G$, where $G$ is a finite Abelian group containing elements of order 2. Denote by $F$ the subgroup of $G$ generated by all elements of $G$ of order 2. Let $\alpha=(a, \alpha_{G})$  be a topological automorphism of the group $X$. Set    $K={\rm Ker}(I+\alpha_{G})$ and assume that   $F=K$. Let $\xi_1$ and $\xi_2$ be
independent random variables with values in   $X$ and distributions
$\mu_1$ and $\mu_2$ with nonvanishing characteristic functions. It turns out that generally speaking the symmetry of the conditional distribution of the linear form $L_2 =  \xi_1 + \alpha\xi_2$ given $L_1 =  \xi_1 +  \xi_2$ does not imply that $\mu_j$ are shifts of convolutions of Gaussian distributions on $\mathbb{R}$ and distributions  supported in $K$. In other words, Theorem \ref{th2} is not true if the group $G$ contains elements of order 2.
Indeed, let $X=\mathbb{R}\times \mathbb{Z}(2)$.  Each topological automorphism   $\alpha$ of the group $X$ is of the form $\alpha(t, g)=(at, g)$, where   $a\ne 0$. It is obvious that if $a\ne -1$, then $K={\rm Ker}(I+\alpha)=\mathbb{Z}(2)$.
The following statement  was proved in  \cite[Lemma 4.1]{F_solenoid}.

{\it  Consider on the group $Y$ the real-valued functions
$$
f_j(s, h) = \begin{cases}\exp\{-\sigma_j s^2\}, &s\in \mathbb{R}, \ h=0,\\ \kappa_j\exp\{-\sigma_j' s^2\}, &s\in \mathbb{R}, \   h=1,
\\
\end{cases}
$$
where     $0<\sigma_j'<\sigma_j$, $0<|\kappa_j|\le\sqrt\frac{\sigma_j'}{\sigma_j}$, $j=1, 2$. Then $f_j(s, h)$ are the characteristic functions of some distributions   $\mu_j$ on the group $X$.}

Let $\alpha$ be a topological automorphism of the group $X$ of the form $\alpha(t, g)=(at, g)$, where $a<0$,   $a\ne -1$. Assume that numbers $\sigma_j, \sigma_j'$ and $\kappa_j$ as above and the  equalities
 \begin{equation}\label{07.03.1}
\sigma_1+a\sigma_2=0, \ \ \sigma'_1+a\sigma'_2=0
\end{equation}
hold. Let
$\xi_1$ and $\xi_2$ be
independent random variables with values in   $X$ and distributions
$\mu_1$ and $\mu_2$. It follows from  (\ref{07.03.1}) that the characteristic functions $\hat\mu_j(s, h)$ satisfy equation (\ref{11.04.1}). By Lemma \ref{lem1},   the conditional distribution of the linear form
 $L_2 =  \xi_1 + \alpha\xi_2$ given $L_1 =  \xi_1 +  \xi_2$
is symmetric. Obviously,  $\mu_j$ can not be represented as  convolutions of Gaussian distributions on $\mathbb{R}$ and distributions supported in  $K$.
We see that, in contrast to Theorem \ref{th2}, the symmetry of the conditional distribution of one linear form of two independent random variables with values in the group $X=\mathbb{R}\times \mathbb{Z}(2)$  given another does not imply that  the distributions of the independent random variables are represented as  convolutions of Gaussian distributions on $\mathbb{R}$ and distributions supported in  $K$.
\end{remark}

\section{ Generalization  of the Heyde  theorem to   groups of the form  $\Sigma_\text{\boldmath $a$}\times G$, where \text{\boldmath $a$}=$(2, 3, 4,\dots)$ and $G$ is a finite Abelian group}

 Consider the sequence \text{\boldmath $a$}= $(2, 3, 4,\dots)$. We will assume that the sequence \text{\boldmath $a$} is fixed throughout this section. Denote by  $\Sigma_\text{\boldmath $a$}$  the corresponding ${\text{\boldmath $a$}}$-adic solenoid (\!\!\cite[(10.12)]{Hewitt-Ross}). The group $\Sigma_\text{\boldmath $a$}$  is compact,  connected  and has dimension 1 (\!\!\cite[(10.13), (24.28)]{Hewitt-Ross}). The character group of the group   $\Sigma_\text{\boldmath $a$}$ is topologically isomorphic to the additive group of rational numbers $\mathbb{Q}$  considering in the discrete topology (\!\!\cite[(25.3), (25.4)]{Hewitt-Ross}). We note that  the ${\text{\boldmath $a$}}$-adic solenoid $\Sigma_\text{\boldmath $a$}$ is a torsion-free group, and for each positive integer       $n$ the ${\text{\boldmath $a$}}$-adic solenoid $\Sigma_\text{\boldmath $a$}$ is a group with unique division by $n$. Hence, multiplication by any rational number is defined in the group $\Sigma_\text{\boldmath $a$}$.  Each topological automorphism of the group $\Sigma_\text{\boldmath $a$}$ is multiplication by a nonzero rational number, and the converse also holds.

It follows from Definition \ref{d1} that the characteristic function of a Gaussian distribution   $\gamma$ on   the  group
 $\Sigma_\text{\boldmath $a$}$ is of the form
$$
\hat\gamma(r)=(z, r)\exp\{-\sigma r^2\}, \ \ r\in \mathbb{Q},
$$
where $z\in  \Sigma_\text{\boldmath $a$}$, $\sigma \ge 0$.

Let $X=\Sigma_\text{\boldmath $a$}\times G$, where   $G$ is a finite Abelian group. Denote by $x=(z, g)$, where $z\in \Sigma_\text{\boldmath $a$}$, $g\in G$, elements of the group    $X$.  Let $H$ be the character group of the group    $G$. The   group  $Y$ is topologically isomorphic to the group   $\mathbb{Q}\times H$. Denote by $y=(r, h)$, where $r\in \mathbb{Q}$, $h\in H$, elements of the group $Y$. Let  $\alpha$ be a topological automorphism of the group
$X$.  Since $\Sigma_\text{\boldmath $a$}$ is the  connected component of zero of $X$ and $G$ is the torsion subgroup of $X$, we have   $\alpha(\Sigma_\text{\boldmath $a$})=\Sigma_\text{\boldmath $a$}$ and  $\alpha(G)=G$. Therefore,   $\alpha$ acts on elements of $X$ as follows $\alpha(z, g)=(\alpha_{{\Sigma_\text{\boldmath $a$}}} z, \alpha_{G} g)$, where $z\in \Sigma_\text{\boldmath $a$}$, $g\in G$. In so doing,  $\alpha_{{\Sigma_\text{\boldmath $a$}}}$ is multiplication by a nonzero rational number $a$.   We will identify $\alpha_{{\Sigma_\text{\boldmath $a$}}}$ with a rational number  $a$.   Thus, $\alpha(z, g)=(a z, \alpha_{G} g)$, and we will write    $\alpha$   in the form $\alpha=(a, \alpha_{G})$. The adjoint  automorphism   $\tilde\alpha_{{\Sigma_\text{\boldmath $a$}}}$ of the group $\mathbb{Q}$ we shall also   identify with $a$.

The following  generalization of   Theorem \ref{th1} holds.

\begin{theorem}\label{th6} Let $X=\Sigma_\text{\boldmath $a$}\times G$,   where \text{\boldmath $a$}= $(2, 3, 4,\dots)$, and $G$ is a finite Abelian group containing no elements of order  $2$.
Let  $\alpha=(a, \alpha_{G})$ be a topological automorphism of the group
  $X$. Set $K={\rm Ker}(I+\alpha_{G})$.
Let $\xi_1$ and $\xi_2$ be
independent random variables with values in   $X$ and distributions
$\mu_1$ and $\mu_2$ with nonvanishing characteristic functions. Assume that the conditional  distribution of the linear form $L_2 = \xi_1 + \alpha\xi_2$ given $L_1 = \xi_1 + \xi_2$ is symmetric.
If  $a \ne-1$, then $\mu_j=\gamma_j*\omega*E_{x_j}$, where $\gamma_j$ is a Gaussian distribution on   $\Sigma_\text{\boldmath $a$}$, $\omega$ is a distribution supported in    $K$, $x_j\in X$, $j=1, 2$.  If  $a=-1$, then $\mu_j=\omega*E_{x_j}$, where  $\omega$ is a distribution supported in $\Sigma_\text{\boldmath $a$}\times K$, $x_j\in X$, $j=1, 2$.
\end{theorem}
The proof of Theorem \ref{th6} is based on the following lemmas.
 \begin{lemma}\label{lem12}  {\rm(\!\!\cite[Theorem 3]{FeTVP1})}    Let  $X$ be a locally compact Abelian group containing no elements of order  $2$, and let  $\alpha$ be a topological automorphism of the group  $X$ satisfying the condition ${\rm Ker}(I+\alpha)=\{0\}$. Let $\xi_1$ and $\xi_2$ be
independent random variables with values in   $X$ and distributions
$\mu_1$ and $\mu_2$ with nonvanishing characteristic functions. Assume that the conditional  distribution of the linear form $L_2 = \xi_1 + \alpha\xi_2$ given $L_1 = \xi_1 + \xi_2$ is symmetric. Then  $\mu_j$  are Gaussian distributions.
\end{lemma}
\begin{lemma}\label{lem4}  Let $X=\Sigma_\text{\boldmath $a$}\times K$,  where \text{\boldmath $a$}= $(2, 3, 4,\dots)$, and $K$ is a locally compact Abelian group. Denote by $(z, k)$, where $z\in \Sigma_\text{\boldmath $a$}$, $k\in K$, elements of the group $X$, and by $L$ the character group of the group $K$. Assume that $L^{(2)}=L$.
Let  $\alpha$ be a topological automorphism of the group
$X$ of the form $\alpha(z, k)=(az, -k)$, where $a\ne -1$.
Let $\xi_1$ and $\xi_2$ be
independent random variables with values in   $X$ and distributions
$\mu_1$ and $\mu_2$ with nonvanishing characteristic functions. Assume that the conditional  distribution of the linear form $L_2 = \xi_1 + \alpha\xi_2$ given $L_1 = \xi_1 + \xi_2$ is symmetric.
Then $\mu_j=\gamma_j*\omega$, where $\gamma_j$ is a Gaussian distribution on  $\Sigma_\text{\boldmath $a$}$, and $\omega$ is a distribution supported in    $K$,   $j=1, 2$.
\end{lemma}
{\rmfamily   {\bfseries   {\itshape Proof  }}}     The group $Y$ is topologically isomorphic to the group   $\mathbb{Q}\times L$. Denote by $y=(r, l)$,  where $r\in \mathbb{Q}$, $l\in L$, elements of the group $Y$. By Lemma \ref{lem1}, the characteristic functions   $\hat\mu_j(r, l)$   satisfy equation (\ref{11.04.1})  which takes the form
$$
\hat\mu_1(r_1+r_2, l_1+l_2)\hat\mu_2(r_1+a r_2, l_1- l_2)$$\begin{equation}\label{03.05.1}=
\hat\mu_1(r_1-r_2, l_1-l_2)\hat\mu_2(r_1-a r_2, l_1+l_2), \ \  r_j\in \mathbb{Q}, \ \ l_j\in L.
\end{equation}
Substitute  $l_1=l_2=0$ in (\ref{03.05.1}). Taking into account that multiplication by a nonzero rational number is a topological automorphism of the group $\Sigma_\text{\boldmath $a$}$  and  Lemma \ref{lem1}, it follows from the   obtained equation and Lemma \ref{lem12}, applying to the group $\Sigma_\text{\boldmath $a$}$ and the topological automorphism $a$, that
$$
\hat\mu_j(r, 0)=(z_j, r)\exp\{-\sigma_jr^2\}, \ \ r\in \mathbb{Q},
$$  where $\sigma_j\ge 0$, $z_j\in \Sigma_\text{\boldmath $a$}$, $j=1, 2$. Since $z_1+az_2=0$, we can replace the distributions $\mu_j$ by  their shifts $\lambda_j=\mu_j*E_{-z_j}$ and suppose from the   beginning,   without loss of generality,      that  $z_1=z_2=0$. We also note that  $\sigma_1+a\sigma_2=0$. This implies that either  $\sigma_1=\sigma_2=0$ or $\sigma_1>0$ and $\sigma_2>0$. The case when $\sigma_1=\sigma_2=0$,    is considered in the same way as in the proof of Lemma  \ref{lem10}. Thus, we will assume that
 \begin{equation}\label{03.05.2}
\hat\mu_j(r, 0)=\exp\{-\sigma_jr^2\}, \ \ r\in \mathbb{Q},
\end{equation}
where $\sigma_j> 0$, $j=1, 2$. Embed in a natural way $(r, l)\rightarrow (r, l)$ the group $\mathbb{Q}\times L$ into the group   $\mathbb{R}\times L$.

Note now that the inequality
\begin{equation}\label{03.05.3}
|\hat\mu(u)-\hat\mu(v)|\le \sqrt 2(1-{\rm Re}\
\hat\mu(u-v))^\frac{1}{2}, \ \ u, v\in Y,
\end{equation}
holds for an arbitrary characteristic function $\hat\mu(y)$ on the group $Y$.   
It follows from  (\ref{03.05.2}) and (\ref{03.05.3}) that the characteristic functions  $\hat\mu_j(r, l)$, $j=1, 2$, are uniformly continuous on the subgroup  $\mathbb{Q}\times L$ in the topology induced on  $\mathbb{Q}\times L$ by the topology of the group   $\mathbb{R}\times L$. Hence, the characteristic functions $\hat\mu_j(r, l)$ can be extended by continuity from the subgroup   $\mathbb{Q}\times L$ to the group $\mathbb{R}\times L$. We retain the notation $\hat\mu_j(s, l)$, where $s\in \mathbb{R}$, $l\in L$, for the continued functions.  It follows from  (\ref{03.05.1}) that the characteristic functions  $\hat\mu_j(s, l)$ satisfy equation (\ref{14.04.2}).   In view of Lemma \ref{lem1}, by Lemma \ref{lem10}, the representation (\ref{14.04.11}) holds. Hence, $$\hat\mu_j(r, l)=\exp\{-\sigma_j r^2\}\hat\omega(l), \ \ r\in \mathbb{Q}, \ \ l\in L, \ \ j=1, 2.$$ The statement of the lemma follows from this. $\hfill\Box$

\medskip

{\rmfamily   {\bfseries   {\itshape Proof  of Theorem \ref{th6}}}}    By Lemma \ref{lem1}, the characteristic functions   $\hat\mu_j(r, h)$   satisfy equation (\ref{11.04.1})  which takes the form
$$
\hat\mu_1(r_1+r_2, h_1+h_2)\hat\mu_2(r_1+a r_2, h_1+\tilde\alpha_{G} h_2)$$\begin{equation}\label{10.05.1}=
\hat\mu_1(r_1-r_2, h_1-h_2)\hat\mu_2(r_1-ar_2, h_1-\tilde\alpha_{G} h_2), \ \  r_j\in \mathbb{Q}, \ \ h_j\in H.
\end{equation}
The proof of the theorem boils down to solving of  equation (\ref{10.05.1}) and repeats the proof of Theorem  \ref{th2}. We only   change   $\mathbb{R}$ for $\Sigma_\text{\boldmath $a$}$ and   use Lemma  \ref{lem4} instead of Lemma \ref{lem10}.

Note that arguing in the same way as in Remark \ref{nr1}, it is easy to verify that   Theorem \ref{th6} can not be strengthened   by narrowing the class of distributions which are characterized by the symmetry of the conditional distribution of one linear form given another.
$\hfill\Box$

Note that if in Theorem \ref{th6} $a\ne -1$, then ${\rm Ker}(I+\alpha)=K$, and if $a=-1$, then ${\rm Ker}(I+\alpha)=\Sigma_\text{\boldmath $a$}\times K$. Thus, Theorem \ref{th6} implies the following statement.
\begin{corollary}\label{co2} Let $X=\Sigma_\text{\boldmath $a$}\times G$,   where \text{\boldmath $a$}= $(2, 3, 4,\dots)$, and $G$ is a finite Abelian group containing no elements of order  $2$.
Let  $\alpha$ be a topological automorphism of the group
  $X$.
Let $\xi_1$ and $\xi_2$ be
independent random variables with values in   $X$ and distributions
$\mu_1$ and $\mu_2$ with nonvanishing characteristic functions. If the conditional  distribution of the linear form $L_2 = \xi_1 + \alpha\xi_2$ given $L_1 = \xi_1 + \xi_2$ is symmetric, then $\mu_j=\gamma_j*\omega$, where $\gamma_j$ is a Gaussian distribution on   $X$, $\omega$ is a distribution supported in     ${\rm Ker}(I+\alpha)$,   $j=1, 2$.
\end{corollary}
In conclusion,   we formulate the following  problem.

\bigskip

\noindent{\bf Problem.}  {\it Let  $X$ be a   locally compact Abelian group containing no elements of order  $2$.
Let  $\alpha$ be a topological automorphism of the group
  $X$. Set $K={\rm Ker}(I+\alpha)$. Let $\xi_1$ and $\xi_2$ be independent random variables with   values in   $X$ and distributions
$\mu_1$ and $\mu_2$ with nonvanishing characteristic functions.
Describe  the class $\Upsilon$ of  such groups  $X$ for which the symmetry of the conditional distribution of the linear form    $L_2 = \xi_1 + \alpha\xi_2$ given $L_1 = \xi_1 + \xi_2$ implies that $\mu_j=\gamma_j*\omega$, where $\gamma_j$ is a Gaussian distribution on   $X$, $\omega$ is a distribution supported in    $K$,  $j=1, 2$.}

\bigskip

 Recall that each  Gaussian distribution on   a locally compact Abelian group $X$ is concentrated in a coset of the connected component of zero of $X$  (\!\!\cite[Chapter IV, \S 6, Remark 2]{Pa}). Taking this into account,     Theorem   \ref{th1} implies that  finite Abelian groups containing no elements of order  $2$ belong to $\Upsilon$.  It follows from  Corollaries  \ref{co1} and \ref{co2} that  the   groups of the form $\mathbb{R}\times G$ and $\Sigma_\text{\boldmath $a$}\times G$, where  $G$ is a finite Abelian group containing no elements of order  $2$, and $\Sigma_\text{\boldmath $a$}$ is the  \text{\boldmath $a$}-adic solenoid with \text{\boldmath $a$}=$(2, 3, 4,\dots)$ also belong to $\Upsilon$. As proven in \cite[Theorem 2.1]{F_solenoid}, $\Upsilon$ contains all  \text{\boldmath $a$}-adic solenoids with no elements of order  $2$ as well.
 
\bigskip

%\newpage

\bigskip

\noindent B. Verkin Institute for Low Temperature Physics and Engineering\\
of the National Academy of Sciences of Ukraine\\
47, Nauky ave, Kharkiv, 61103, Ukraine

\medskip

\noindent Department of Mathematics  
University of Toronto \\
40 St. George Street
Toronto, ON,  M5S 2E4
Canada 

\medskip

\noindent e-mail:    gennadiy\_f@yahoo.co.uk

\end{document}